\documentclass[12pt]{article}
\usepackage[margin=.75in]{geometry}            

\usepackage[intlimits]{amsmath}
\usepackage{amsfonts,amssymb}
\DeclareSymbolFontAlphabet{\mathbb}{AMSb}
\usepackage{float}
\usepackage[bf]{caption}       
\setcaptionmargin{0.5in}
\usepackage{fancyhdr}
\usepackage{fancybox}
\usepackage{ifthen}
\usepackage{url}
\usepackage{lscape,afterpage}
\usepackage{xspace}
\usepackage{setspace}

\usepackage[colorlinks=true, linkcolor=blue, citecolor=blue, urlcolor=blue]{hyperref}
\usepackage[shortlabels, inline]{enumitem}
\usepackage{amsthm}
\usepackage{mathtools}
\usepackage{color}
\usepackage[color = lime]{todonotes}
\usepackage{tabularx}
\usepackage{diagbox}
\usepackage{tikz-cd}
\usepackage{graphicx}

\theoremstyle{plain}
\newtheorem{thm}{Theorem}[section]
\newtheorem{conj}[thm]{Conjecture}
\newtheorem{lemma}[thm]{Lemma}
\newtheorem{prop}[thm]{Proposition}

\theoremstyle{definition}
\newtheorem{obs}[thm]{Remark}
\newtheorem{mydef}[thm]{Definition}
\newtheorem{ex}[thm]{Example}

\newcommand{\A}{\mathbb{A}}
\newcommand{\F}{\mathbb{F}}
\newcommand{\PP}{\mathbb{P}}
\newcommand{\Q}{\mathbb{Q}}
\newcommand{\C}{\mathbb{C}}
\newcommand{\Z}{\mathbb{Z}}
\newcommand{\E}{\mathcal{E}}

\newcommand{\OO}{\mathcal{O}}
\newcommand{\rk}{\textnormal{rank}}
\newcommand{\defeq}{\mathrel{\mathop:}=}
\newcommand{\dsty}{\displaystyle}
\newcommand{\Sht}{\textnormal{Sht}}

\newcommand{\GL}{\textnormal{GL}}
\newcommand{\tr}{\textnormal{tr}}
\newcommand{\overbar}[1]{\mkern 1.5mu\overline{\mkern-3mu#1\mkern-0.5mu}\mkern 1.5mu}

\title{Moduli Spaces of Shtukas over the Projective Line}

\author{Mar\'ia In\'es de Frutos--Fern\'andez}

\date{}

\begin{document}

\maketitle

\begin{abstract}
We provide explicit equations for moduli spaces of Drinfeld shtukas over the projective line with $\Gamma(N)$, $\Gamma_1(N)$ and $\Gamma_0(N)$ level structures, where $N$ is an effective divisor on $\PP^1$. If the degree of $N$ is high enough, these moduli spaces are relative surfaces. We study some invariants of the moduli space of shtukas with $\Gamma_0(N)$ level structure for several degree $4$ divisors on $\PP^1$.
\end{abstract}

\section{Introduction}
Shtukas were introduced by Drinfeld in the seventies, in his series of papers containing the proof of the Langlands correspondence for GL$_2$ over function fields \cite{Dri_pf, Dri_Lan}, based on the study of the cohomology of moduli spaces of rank $2$ shtukas.

Let $X$ be a smooth, projective, geometrically irreducible curve over $\F_q$ with  function field $K$. Denote $\overbar{X} := X \times_{\F_q} \overline{\F}_q$, where $\overline{\F}_q$ is an algebraic closure of $\F_q$.
Given any scheme $S$ over $\F_q$ and any vector bundle $\E$ over $X \times S$, we denote by 
$\mathcal{E}^\sigma$ the pullback $(\text{Id}_X \times \text{Frob}_S)^*\mathcal{E}$, where $\text{Frob}_S : S \to S$ is the
Frobenius morphism which is the identity on points, and $t \mapsto t^q$ on functions.

\begin{mydef}
	Let $S$ be a scheme over $\mathbb{F}_q$. A \textit{(Drinfeld) shtuka of rank $r$ over $S$} is a tuple $\tilde{\mathcal{E}} = ( \mathcal{E}, P, Q, \mathcal{E} \hookrightarrow \mathcal{E}' \hookleftarrow \mathcal{E}^\sigma)$, consisting of
	\begin{itemize}[--]
		\item a rank $r$ vector bundle $\mathcal{E}$ over $X \times S$,
		\item two morphisms $P, Q : S \to X$, called the \textit{pole} and the \textit{zero} of the shtuka $\tilde{\mathcal{E}}$,
		\item a modification consisting of two injections
		$\mathcal{E} \hookrightarrow \mathcal{E}' \hookleftarrow \mathcal{E}^\sigma,$ where $\mathcal{E}'$ is a rank $r$ vector bundle over $X \times S$ such that $\mathcal{E}'/\mathcal{E}$ and $\mathcal{E}'/\mathcal{E}^\sigma$ are supported on the graphs of $P$ and $Q$, respectively, and they are invertible on their support.
	\end{itemize}
	Two shtukas  $\tilde{\mathcal{E}} = ( \mathcal{E}, P, Q, \mathcal{E} \hookrightarrow \mathcal{E}' \hookleftarrow \mathcal{E}^\sigma)$ and $\tilde{\mathcal{F}} = ( \mathcal{F}, P', Q', \mathcal{F} \hookrightarrow \mathcal{F}' \hookleftarrow \mathcal{F}^\sigma)$ of rank $r$ are \textit{isomorphic} when $P = P'$, $Q = Q'$, and there exist isomorphisms $\alpha : \E \to \mathcal{F}$ and $\beta : \E' \to \mathcal{F}'$ such that the following diagram commutes:
	\[
	\begin{tikzcd}[column sep = 1.5cm]
	\E \arrow[r, hook] \arrow{d}{\alpha} & \arrow{d}{\beta} \E'& \arrow[hook']{l} \arrow{d}{\alpha^\sigma} \E^\sigma \\   
	\mathcal{F} \arrow[r, hook]  & \mathcal{F}' & \arrow[hook']{l} \mathcal{F}^\sigma.
	\end{tikzcd}
	\]
	
\end{mydef}

\begin{mydef}
	The \textit{stack} Sht$^r \defeq$ Sht$^r_X$ \textit{of shtukas of rank $r$} associates to a scheme $S$ over $\mathbb{F}_q$ the groupoid Sht$^r(S)$ whose objects are shtukas of rank $r$ over $S$.
\end{mydef}

Drinfeld proved that the stack $\Sht^r$ is an algebraic Deligne--Mumford
stack, and the map $ (P, Q) : \Sht^r \to X \times X $
sending a shtuka to its pole and zero is smooth of relative dimension $2r - 2$. Moreover, $\Sht^r$ is locally of finite type (see \cite{Dri_Lan} or \cite[Theorem 2.2]{Tua}).

Denote by $\Sht^{r, d}$ the substack of $\Sht^r$ consisting of shtukas for which the vector bundle $\E$ has degree $d$. Then each $\Sht^{r, d}$ is a connected component of $\Sht^{r}$, and
\[ \Sht^{r} = \bigsqcup_{d \in \Z} \Sht^{r, d}.\]

Let $S$ be an $\F_q$--scheme. Let $N$ be a finite subscheme of $X$ and let $\tilde{\E}$ be a shtuka of rank $r$ over $S$ such that the graphs of its pole $ P : S \to X$ and zero $Q : S \to X$ do not intersect $N \times S$. Then the data of the shtuka provides a restricted isomorphism $ \psi : \E_{| _{N \times S}} \to \E^\sigma_{| _{N \times S}}.$

\begin{mydef}
	A \textit{level structure of $\tilde{\E}$ on $N$} is an isomorphism
	$ \alpha:  \OO^r_{N \times S} \to \E_{| _{N \times S}} $
	such that the following diagram commutes:
	\[
	\begin{tikzcd}[column sep = 3cm]
	\OO^r_{N \times S} \arrow{r}{=} \arrow{d}{\alpha} & \arrow{d}{\alpha^\sigma} \OO^{r,\sigma}_{N \times S} \\   
	\E_{| _{N \times S}} \arrow{r}{\psi}  & \E^\sigma_{| _{N \times S}}
	\end{tikzcd}
	\]
\end{mydef}

We denote by Sht$^r(N)$ the \textit{stack of
Drinfeld shtukas of rank $r$ with level structure on $N$}, which associates to a scheme $S$ over $\mathbb{F}_q$ the groupoid Sht$^r(N)(S)$ whose objects are pairs $(\tilde{\E}, \alpha)$, where $\tilde{\E}$ is a rank $r$ shtuka over $S$ and $\alpha$ is a level structure of $\tilde{\E}$ on $N$, with the evident notion of isomorphism.

Drinfeld proved that the forgetful morphism $\Sht^r(N) \to \Sht^r \times_{X^2} (X - N)^2$	is representable, finite, \'etale and Galois with Galois group $\textnormal{GL}_ r(\OO_N)$ (see \cite[Section 3]{Dri_var}).

Now we restrict our attention to the case where the shtuka $\tilde{\E}$ has rank $2$. Then we often use the phrase ``\textit{$\Gamma(N)$ level  structure of $\tilde{\E}$}'' to refer to a level structure on $N$, and use Sht$^2(\Gamma(N))$ to denote Sht$^2(N)$. We introduce two new notions of level structure:
\begin{enumerate}[\normalfont(1)]	
	\item A \textit{$\Gamma_1(N)$ level  structure of $\tilde{\E}$} is a section $v \in H^0(N\times S, \E_{|_{N \times S}})$ which generates a locally free $\OO_{N \times S}$--module of rank $1$ and such that $v^\sigma = \psi(v) $.
	\item A \textit{$\Gamma_0(N)$ level  structure of $\tilde{\E}$} is a rank $1$ locally free $\OO_{N \times S}$--submodule $\mathcal{L} \subset \E_{|_{N \times S}}$ such that $\mathcal{L}^\sigma = \psi(\mathcal{L}) $.
\end{enumerate}

We denote by Sht$^2(\Gamma_1(N))$ the \textit{stack of
Drinfeld shtukas of rank $2$ with $\Gamma_1(N)$ level  structure}, which associates to a scheme $S$ over $\mathbb{F}_q$ the groupoid Sht$^2(\Gamma_1(N))(S)$ whose objects are pairs $(\tilde{\E}, v)$, where $\tilde{\E}$ is a shtuka of rank $r$ over $S$ and $v$ is a $\Gamma_1(N)$ level structure. We define Sht$^2(\Gamma_0(N))$ analogously.

\begin{obs}
	For some purposes, the study of the stacks $\Sht^r$ or $\Sht^r(\Gamma)$ for some level structure $\Gamma$ is not enough, and one must instead use the compactified stacks $\overline{\Sht}^{\:r}$ or $\overline{\Sht}^{\:r}(\Gamma)$. To create these compactifications, the notion of generalized shtuka is introduced. For $\Gamma = \Gamma(N)$, the compactification was studied by Drinfeld in the rank $2$ case \cite{Dri_coh}, and generalized to arbitrary rank by Lafforgue \cite{Laff}. It played a key role in the proof of the Langlands correspondence for GL$_n$ over function fields.
\end{obs}

%MAIN RESULTS

The main goal of this paper is to derive explicit equations for moduli spaces of shtukas of rank $2$ over the base curve $X = \PP^1_{\F_q}$, with the kinds of level structures mentioned above. 
Using these equations, we conclude the results summarized in theorems \ref{thm-main-1} and \ref{thm-main-2}.

\begin{thm}\label{thm-main-1}
	Let $N$ be an effective divisor on $\PP^1$. Denote by $\Delta$ the diagonal subscheme of $(\PP^1 \setminus N) \times (\PP^1 \setminus N)$.
	
	\begin{enumerate}[\normalfont(1)]
		\item If $\deg N \geq 1$, then $\Sht^{2, \tr}(\Gamma(N))$ is birational to a surface over $((\PP^1 \setminus N) \times (\PP^1 \setminus N)) \setminus \Delta$.
		\item If $\deg N \geq 2$, then $\Sht^{2, \tr}(\Gamma_1(N))$ is birational to a surface over $((\PP^1 \setminus N) \times (\PP^1 \setminus N)) \setminus \Delta$.
		\item If $\deg N \geq 3$, then $\Sht^{2, \tr}(\Gamma_0(N))$ is birational to the stacky quotient of a surface over $((\PP^1 \setminus N) \times (\PP^1 \setminus N)) \setminus \Delta$ by the trivial action of the multiplicative group $\mathbb{G}_m$.
	\end{enumerate}
\end{thm}

\begin{thm}\label{thm-main-2}
	Let $N$ be an effective divisor on $\PP^1$ supported on points of degree $1$. Then up to birational equivalence, we have that:
	
	\begin{enumerate}[\normalfont(1)]
		\item If $\deg N = 1$, then $\Sht^{2, \tr}(\Gamma(N))$ is a rational surface.
		\item If $\deg N = 2$, then $\Sht^{2, \tr}(\Gamma_1(N))$ is a rational surface.
		\item If $\deg N = 3$, then $\Sht^{2, \tr}(\Gamma_0(N))$ is the quotient of a rational surface by $\mathbb{G}_m$.
		\item If $\deg N = 4$, then $\Sht^{2, \tr}(\Gamma_0(N))$ is the quotient of an elliptic surface by $\mathbb{G}_m$.
	\end{enumerate}
\end{thm}

%CONTENTS
The contents of this paper are as follows. In Section \ref{sec-surf} we review some basic background on algebraic surfaces. In Section  \ref{sec-sht} we study moduli spaces of shtukas of rank $2$ over the projective line; we describe the stack of shtukas with no level structure and with $\Gamma(N)$, $\Gamma_1(N)$, and $\Gamma_0(N)$ level structures, giving explicit equations for these spaces.
Section \ref{sec-deg-4} is dedicated to the study of some arithmetic invariants of $\Sht^{2, \tr}(\Gamma_0(N))$ for $N$ a degree $4$ effective divisor on $\PP^1$.
Finally, in Section \ref{sec-fw} we describe some future and related work in this topic, including a modularity conjecture for elliptic curves over funcion fields involving the moduli stack of shtukas with $\Gamma_0(N)$ level structure.

\section{Background on Surfaces}\label{sec-surf}
We review some definitions and results about surfaces which will be used in sections \ref{sec-sht} and \ref{sec-deg-4},
mainly following \cite{ScSh}. In this section $X$ denotes a smooth projective curve over a field $k$.

\begin{mydef}
	A \textit{rational surface} is a surface which is birational to the projective plane.
	An \textit{elliptic surface} $S$ over $X$ is a smooth projective surface $S$ with an elliptic fibration over $X$, that is, with a surjective morphism $f : S \to X$ whose base change to the algebraic closure is such that almost all fibers are smooth curves of genus $1$ and no fiber contains a smooth rational curve of self--intersection $-1$.
\end{mydef}

A \textit{section} of an elliptic surface $f: S \to X$ is a morphism $\pi: X \to S$
	such that $f \circ \pi = \text{id}_X$.
If $S$ has a section, then we can regard its generic fiber $E$ as an elliptic curve over the field $k(X)$, which allows us to work with a Weierstrass model.
In this paper we only consider elliptic surfaces $S$ which admit a section and have at least one singular fiber.

The classification of singular fibers of elliptic surfaces was carried out by Kodaira \cite{Kod, Kod2, Kod3} in the case where $k$ is the complex numbers, and generalized to other fields by N\'eron \cite{Ner}. The imperfect residue field case was studied by Szydlo \cite{Szy}.

The N\'eron--Severi group $\textnormal{NS}(S)$ of a surface $S$ is the quotient of its divisor group under algebraic equivalence. The rank of the N\'eron--Severi group is called the \textit{Picard number} of $S$, denoted $\rho(S)$.

If $S$ is an elliptic surface over $X$, this group is closely related to the Mordell--Weil group $E(K)$ of its generic fiber $E$, where $K := k(X)$. Denote by $T$ the \textnormal{trivial lattice of $S$}, that is, the subgroup of the N\'eron--Severi group of $S$ generated by the zero section and the fiber components. Then
$E(K)$ is isomorphic to $\textnormal{NS}(S)/T$ \cite[Theorem 6.3]{ScSh}.
In particular,
\[ \rho(S) = \rk(T) + \textnormal{rank}_{MW}(E) . \]

The Mordell--Weil rank of $E$ is in general hard to compute, but the rank of the trivial lattice can be computed from local information: denote by $m_\mathit{v}$ the number of irreducible components of each fiber $F_\mathit{v}$ of $S$. Then
	\[ \rk(T) = 2 + \sum_{\mathit{v} \in \overbar{X}} (m_\mathit{v} - 1).\]

\begin{mydef}
	The \textit{arithmetic genus} of a projective surface $S$ over $k$ is defined as
	\[ p_a(S) := \sum_{i=0}^{2} (-1)^{i}\dim_k H^i(S, \OO_S).\]
\end{mydef}

If $S$ is a complete intersection inside a product of projective spaces such that $S$ is smooth or has at most rational double points\footnote{An isolated singularity of $S$ is called a \textit{rational double point} or \textit{du Val singularity} if there exists a resolution of singularities $\psi : S' \to S$
	such that $K_{S'} = \psi^*K_S$, where $K_S$ denotes the canonical divisor of $S$.}, then its arithmetic genus depends only on the degrees of the polynomials describing it. In particular:
\begin{lemma}\label{lemm-arith-gen}
	Let $S \subset \PP:= \PP^1 \times \PP^1 \times \PP^1$ be a hypersurface given by a  multi--homogeneous polynomial of multi--degree $d_1, d_2, d_3$, with $d_i \geq 2$. 
	Assume that $S$ is smooth or has at most rational double points.
	Then the arithmetic genus of $S$ is $(d_1-1)(d_2-1)(d_3-1) + 1$.
\end{lemma}
\begin{proof}
	Consider the cohomology long exact sequence of the closed subscheme exact sequence 
	\[ 0 \to \mathcal{O}_\PP(-S) = \mathcal{O}_\PP(-d_1, -d_2, -d_3) \to \mathcal{O}_\PP \to \iota_*\mathcal{O}_S \to 0 \]
	corresponding to the inclusion $\iota : S \to \PP$.
	From K\"unneth's formula, 
	\[ H^n(\PP, \mathcal{O}_\PP) = \bigoplus_{i + j + k = n} H^i(\PP^1, \mathcal{O}_{\PP^1}) \otimes H^j(\PP^1, \mathcal{O}_{\PP^1}) \otimes H^k(\PP^1, \mathcal{O}_{\PP^1}), \]
	so that
	\[ h^n(\PP, \mathcal{O}_\PP) := \dim H^n(\PP, \mathcal{O}_\PP) = \begin{cases}
	1, \text{ if } n=0,\\
	0, \text{ if } n \geq 1.
	\end{cases}\]
	Analogously,
	\begin{align*}
	H^n(\PP, \mathcal{O}_\PP(-d_1, -d_2, -d_3)) 
	 = \bigoplus_{i + j + k = n} H^i(\PP^1, \mathcal{O}_{\PP^1}(-d_1)) \otimes H^j(\PP^1, \mathcal{O}_{\PP^1}(-d_2)) \otimes H^k(\PP^1, \mathcal{O}_{\PP^1}(-d_3)). 
	\end{align*}
	Hence $H^n(\PP, \mathcal{O}_\PP(-d_1, -d_2, -d_3)) = 0$ for $n=0, 1, 2$, since at least one of $i, j, k$ must be zero in each summand of the formula. However,
	\begin{align*}
	H^3(\PP, \mathcal{O}_\PP(-d_1, -d_2, -d_3)) 
	&= H^1(\PP^1, \mathcal{O}_{\PP^1}(-d_1)) \otimes H^1(\PP^1, \mathcal{O}_{\PP^1}(-d_2)) \otimes H^1(\PP^1, \mathcal{O}_{\PP^1}(-d_3)) =\\
	&= H^0(\PP^1, \mathcal{O}_{\PP^1}(d_1-2)) \otimes H^0(\PP^1, \mathcal{O}_{\PP^1}(d_2-2)) \otimes H^0(\PP^1, \mathcal{O}_{\PP^1}(d_3-2)),
	\end{align*} 
	where the second equality is obtained by Serre duality. Since each $d_i$ is at least $2$, it follows that
	\[h^3(\PP, \mathcal{O}_\PP(-d_1, -d_2, -d_3)) = (d_1-1)(d_2-1)(d_3-1).\]
	Finally, we use the cohomology exact sequence to obtain $h^i(S,\mathcal{O}_S)$ for $i=0,1,2$, and conclude that the arithmetic genus of $S$ is
	\begin{align*}
	p_a(S) = h^2(S,\mathcal{O}_S) - h^1(S,\mathcal{O}_S) + h^0(S,\mathcal{O}_S) 
	= (d_1-1)(d_2-1)(d_3-1) - 0 + 1.
	\end{align*}
\end{proof}

\begin{mydef}
	Let $S$ be a projective surface over a field $k$, and let $\ell$ be a prime different from the characteristic of $k$. The \textit{Euler number} $e(S)$ of $S$ is
	\[ e(S) := \sum_{i=0}^{2} (-1)^{i}\dim_k H^i_{\text{\'et}}(S, \Q_\ell).\]
\end{mydef}

If $S$ is an elliptic surface, then its Euler number, arithmetic genus, and second Betti number $b_2$ are related by the formulas
	\[ e(S) = 12 p_a(S) > 0 \]
and
\[ b_2(S) = e(S) - 2(1 - b_1(X)). \]
In particular, if $X$ is the projective line, then $b_2(S) = e(S) - 2 = 12p_a(S) - 2$.

We conclude this review by making a few remarks about rational and K3 elliptic surfaces. Rational elliptic surfaces are always fibered over the projective line, and they are characterized by having arithmetic genus $p_a(S) = 1$.

\begin{mydef}
	A smooth surface $S$ is called \textnormal{K3} if its canonical bundle is trivial and 
	$h^1(S, \OO_S) = 0$.
\end{mydef}

If $S$ is an elliptic surface over the projective line, then it is K3 if and only if $p_a(S) = 2$. K3 elliptic surfaces are the only kind that can admit two distinct elliptic fibrations with section which are not of product type.

\begin{ex}
	Let $E, E'$ be two elliptic curves over a field $k$. The Kummer surface $S := \text{Km}(E \times E')$ associated to $E \times E'$ is a K3 elliptic surface, and the projections onto $E$ and $E'$ are two elliptic fibrations. The Picard number of $S$ is $ \rho(S) := 18 + \rk \text{ Hom}(E, E')$.
	If $k = \overline{\F}_p$, then the Tate conjecture, known for K3 elliptic surfaces, implies that the Picard number is even. Since $b_2(S)=22$, it follows that the possible values for $\rho(S)$ are $18, 20$, or $22$.
	
\end{ex}

\section{Moduli Spaces of Shtukas of Rank 2 over $\mathbb{P}^1$}\label{sec-sht}

\subsection{The Moduli Stack of Shtukas of Rank 2 over $\mathbb{P}^1$}\label{sec-sht-p1}

Recall that we can describe $\Sht^2 := \Sht^2_{\PP^1}$ as a disjoint union
\[ \Sht^{2} = \bigsqcup_{d \in \Z} \Sht^{2, d},\]
where $\Sht^{2, d}$ is the substack consisting of shtukas for which the vector bundle $\E$ has degree $d$. We will describe the semistable locus of $\Sht^{2, 0}$. Note that the only semistable vector bundle of rank $2$ and degree $0$ over $\PP^1$ is locally on $S$ isomorphic to the trivial one, $\OO \oplus \OO$. This substack of $\Sht^{2, 0}$ is isomorphic to the corresponding substack of $\Sht^{2, d}$ for any even $d \in \Z$. 

Let $\text{Sht}^{2, \tr}$ be the substack of $\text{Sht}^{2, 0}$ whose $S$--points are shtukas with disjoint pole and zero and such that $\E_x$ is trivial on all geometric points $x$ of $S$,
and let $\text{GSht}^{2, \tr}$ be the stack classifying pairs
$\{ (\tilde{\mathcal{E}}, \phi_{\tilde{\mathcal{E}}} :  \OO_{X \times S} \oplus \OO_{X \times S} \xrightarrow{\simeq} \E) \}$ consisting of a shtuka $\tilde{\mathcal{E}}$  with $\E \simeq \OO_{X \times S} \oplus \OO_{X \times S}$, together with a trivialization of $\E$. 
Then $\text{Sht}^{2, \tr}$ is the stacky quotient of $\text{GSht}^{2, \tr}$, which we will now show is represented by a scheme, under the natural action of GL$_2$.

Let $S = \text{Spec}(A)$ be an affine scheme over $\F_q$, and let
$(\tilde{\mathcal{E}} = (\mathcal{E}, P, Q, \mathcal{E} \hookrightarrow \mathcal{E}' \hookleftarrow \mathcal{E}^\sigma), \phi_{\tilde{\mathcal{E}}})$ be an $S$--point of  $\text{GSht}^{2, \tr}$. Note that $P$ and $Q$  are two distinct points in $\PP^1(A)$. Denote $X_A := X \times_{\F_q} \text{Spec}(A)$.

The modification $\mathcal{E} \hookrightarrow \mathcal{E}' \hookleftarrow \mathcal{E}^\sigma$ gives a rational morphism $f: \E \dashrightarrow \E^\sigma$ with a simple pole at $Q$ and well--defined everywhere else, which
induces a morphism $f: \E \to \E^\sigma \otimes \OO(Q)$. Moreover, since the determinant of a shtuka $\tilde{\E}$ is a shtuka $\det \tilde{\E}$ of rank $1$ with the same pole and zero, we must have that div$(\det f) = P - Q$.

Assume that both $P$ and $Q$ lie in the affine patch of $\PP^1$ away from the point at infinity. Then a basis for $H^0(\OO_{X_A}(Q))$ is given by  $\left \{1, \frac{T - P}{T - Q} \right \}$ (if $P = \infty$, then $\frac{T - P}{T - Q}$ gets replaced by $\frac{1}{T - Q}$, and if $Q = \infty$, by $T - P$). Therefore the  morphism $f: \OO_{X_A} \oplus \OO_{X_A} \simeq \E \to \OO_{X_A}(Q) \oplus \OO_{X_A}(Q) \simeq  \E^\sigma \otimes \OO_{X_A}(Q)$ induced by the shtuka $\tilde{\mathcal{E}}$ can be expressed in the given choice of basis by a matrix of the form
\[ M(T) := \begin{bmatrix}
a_0 + a_1 \frac{T - P}{T - Q} & b_0 + b_1 \frac{T - P}{T - Q} \\[6pt]
c_0 + c_1 \frac{T - P}{T - Q} & d_0 + d_1 \frac{T - P}{T - Q}
\end{bmatrix} \]
with $a_0, a_1 \dots, d_0, d_1$ in $A$ and whose determinant is a unit constant multiple of $\frac{T - P}{T - Q}$.

We can identify the subset $\text{GSht}^{2, \tr}_{P, Q}(\text{Spec}(A))$ of $\text{GSht}^{2, \tr}(\text{Spec}(A))$ consisting of pairs in which the shtuka has pole $P$ and zero $Q$ with the set of matrices
\begin{align*}
\mathcal{M}_{A, P, Q} := &\left\{ M(T) := \begin{bmatrix}
a_0 + a_1 \frac{T - P}{T - Q} & b_0 + b_1 \frac{T - P}{T - Q} \\[6pt]
c_0 + c_1 \frac{T - P}{T - Q} & d_0 + d_1 \frac{T - P}{T - Q}
\end{bmatrix} \quad \middle | \quad \begin{matrix}
\\
a_i, b_i, c_i, d_i \in A,\\ \det M(T) = u\frac{T - P}{T - Q}, \\\text{ with } u \in A^*\\
\\
\end{matrix}  \right\} =  \\
=&\left\{ \begin{bmatrix}
a_0 + a_1 \frac{T - P}{T - Q} & b_0 + b_1 \frac{T - P}{T - Q} \\[6pt]
c_0 + c_1 \frac{T - P}{T - Q} & d_0 + d_1 \frac{T - P}{T - Q}
\end{bmatrix} \quad \middle | \quad \begin{matrix}
a_i, b_i, c_i, d_i \in A,\\
a_0d_0 - b_0c_0 = 0,\\
a_0d_1 + a_1d_0 - b_0c_1 - b_1c_0 \in A^*,\\
a_1d_1 - b_1c_1 = 0.\\
\end{matrix} \right\}.
\end{align*}

Hence $\text{GSht}^{2, \tr}$ is represented by the subvariety $V$ of $((\PP^1 \times \PP^1) \setminus \Delta) \times \A^8 $ described by the equations
\[\begin{cases}
a_0d_0 - b_0c_0 = 0,\\
a_0d_1 + a_1d_0 - b_0c_1 - b_1c_0 \in A^*,\\
a_1d_1 - b_1c_1 = 0.\\
\end{cases}
\]

\begin{mydef}
	We say that two matrices $M, N$ in M$_{2\times2}(A(T))$ are \textit{$\sigma$--conjugate} if there exists a matrix $Z \in \text{GL}_2(A)$ such that
	\[ N = Z^\sigma M Z^{-1}.\]
\end{mydef}
GL$_2$ acts on $V$ by  $\sigma$--conjugation, and $\Sht^{2, \tr}$ is the stacky quotient of $V$ under this action.

\subsection{Shtukas with $\mathbf{\Gamma(N)}$ Level Structure}\label{sec-gam}
Let $S = \textnormal{Spec}(A)$ be an affine scheme over $\F_q$, and let $N$ be an effective divisor on $\PP^1$. Let $\tilde{\E}$ be a rank $2$ shtuka over $S$ whose pole $P$ and zero $Q$ do not belong to the support of $N \times S$. Recall that in this situation the data of the shtuka $\tilde{\E}$ induces an isomorphism
$ \psi : \E_{| _{N \times S}} \to \E^\sigma_{| _{N \times S}}, $
and that a $\Gamma(N)$ level structure on $\tilde{\E}$ is an isomorphism
$ \alpha:  \OO^2_{N \times S} \to \E_{| _{N \times S}}$
such that $\alpha^\sigma = \psi \circ \alpha$.

From Section \ref{sec-sht-p1}, we know that after fixing a basis, $\tilde{\E}$ is represented by a matrix in $\mathcal{M}_{A, P, Q}$. Equivalently, if we multiply this matrix by $T - Q$ to cancel the pole at $Q$, we get that $\tilde{\E}$ is represented by a matrix of the form
\[ M(T) := \begin{bmatrix}
a_0 + a_1 T & b_0 + b_1 T \\
c_0 + c_1 T & d_0 + d_1 T
\end{bmatrix}\]
with $\det M(T) = n(T-P)(T-Q)$ for some unit $n \in A^*$.

Under the fixed choice of basis, the level structure can be described by a matrix $\alpha \in \GL_2(\OO_{N \times S})$ such that $\alpha^\sigma = M(N)\alpha$, where $M(N)$ denotes the matrix of the restricted isomorphism $\psi$. In particular, if $N$ is supported on points of degree one, we can make this condition very explicit: at each point $R$ appearing in $N$ with multiplicity $m_R$, we denote by $\pi$ the uniformizer $T-R$ at $R$ and impose the condition
\[\alpha_R^\sigma = M(\pi + R)\alpha_R \mod \pi^{m_R},\]
where $\alpha_R$ denotes the restriction of $\alpha$ to $\OO^2_{(m_R R) \times S} \simeq (\F_q[\pi]/(\pi^{m_R}) \otimes A)^2$. The system of equations thus obtained describes the $\Gamma(N)$ level structure.

\begin{thm}\label{thm-gam-N}
	Let $N$ be a nontrivial effective divisor on $\PP^1$. Then $\Sht^{2, \tr}(\Gamma(N))$ is a surface over $((\PP^1 \setminus N) \times (\PP^1 \setminus N)) \setminus \Delta$, where $\Delta$ is the diagonal subscheme. Moreover,
	\begin{enumerate}[\normalfont(1)]
		\item If $N = (R)$ for a degree $1$ point $R$, then $\Sht^{2, \tr}(\Gamma(N))$ is a rational surface.
	\end{enumerate}
\end{thm}

\begin{proof}
	Let $N$ be a nontrivial effective divisor on $\PP^1_{\F_q}$, and let $R$ be a degree $n$ point in the support of $N$. We will prove that $\Sht^{2, \tr}(\Gamma((R)))$ is a surface over $((\PP^1 \setminus N) \times (\PP^1 \setminus N)) \setminus \Delta$; since the natural map $\Sht^{2, \tr}(\Gamma(N)) \to \Sht^{2, \tr}(\Gamma((R)))$ given by forgetting the extra structure is representable \cite[Proposition 2.3]{Dri_var}, this is enough to conclude the general case.
	
	Let $S = \text{Spec}(A)$ be an affine scheme over $\F_q$ and choose two distinct nonzero points $P, Q \in \PP^1(A)$. An $A$--point of $\Sht^{2, \tr}_{P, Q} (\Gamma((R)))$ is given locally on $S$ by a pair $(M(T), \alpha)$, where $M(T)$ is a matrix of the form
	\[ M(T) := \begin{bmatrix}
	a_0 + a_1 T & b_0 + b_1 T \\
	c_0 + c_1 T & d_0 + d_1 T
	\end{bmatrix}\]
	with determinant $\det M(T) = n(T-P)(T-Q) \textnormal{ for some } n \in A^*$, $\alpha$ is in GL$_2(\OO_{R \times S}) \simeq \text{GL}_2(\F_{q^n} \otimes_{\F_q} A)$, and $\alpha^\sigma = M((R))\alpha$. This pair is considered up to the action of $\text{GL}_2(A)$, which sends $(M(T), \alpha)$ to $(Z^\sigma M(T)Z^{-1}, Z\alpha)$. In particular, if $A = L$ is an algebraically closed field over $\F_q$, then $\OO_{R \times S} \simeq \F_{q^n} \otimes_{\F_q} L \simeq L^n$, and we can identify GL$_2(\OO_{R \times S})$ with $(\GL_2(L))^n$. Suppose that $\alpha$ corresponds to $(\alpha_1, \dots, \alpha_{n})$ under this identification. Then by choosing $Z = \alpha_1^{-1}$, we see that the element of $\Sht^{2, \tr} (\Gamma((R)))$ represented by $(M(T), \alpha)$ has a well--defined representative of the form $(M'(T), \alpha' = (\text{I}_{2\times 2}, \alpha_2, \dots, \alpha_{n}))$. Hence $\Sht^{2, \tr} (\Gamma((R)))$ is a surface.

	To prove claim (1), let $N$ be a degree $1$ effective divisor on $\PP^1$. Possibly after a change of variables, we can assume that $N = (0)$. Fix a ring $A$ over $\F_q$ 
	and two nonzero distinct points $P, Q \in \PP^1(A)$. Then an 
	$A$--point of $\Sht^{2, \tr} (\Gamma((0)))$ is uniquely represented by a pair $(M(T), \alpha = \text{I}_{2 \times 2})$ such that $\text{I}_{2\times 2} = M(0)$, that is, it can be represented by a well--defined matrix of the form
	\[ M(T) := \begin{bmatrix}
	1 + a T & b T \\
	c T & 1 + d T
	\end{bmatrix},\]
	with $\det M(T) = 1 + (a + d)T + (ad - bc)T^2 = n(T-P)(T-Q), n \in A^*$. This  is equivalent to the three conditions
	\[
	\{
	1 = nPQ,\quad
	a + d = -n(P+Q),\quad
	ad - bc = n
	\},
	\]
	which are in turn equivalent to
	\[
	\left\{
	n = \frac{1}{PQ},\quad
	a = -d - \frac{P+Q}{PQ},\quad
	\left(-d - \frac{P+Q}{PQ}\right)d - bc = \frac{1}{PQ}
	\right\}.
	\]
	Hence $\Sht^{2, \tr}(\Gamma((0)))$ is the surface inside $\A^3_{b, c, d}$ cut out by the equation
	\begin{equation}\label{eq-sht-0}
	d^2 + \frac{P+Q}{PQ}d + bc + \frac{1}{PQ} = 0.
	\end{equation}
	
	Denote by $Y$ its projective closure in $\PP^3_{b,c,d,z}$, described by
	\[ d^2 + \frac{P+Q}{PQ}dz + bc + \frac{1}{PQ}z^2 = 0.\]
	
	Using the Jacobian criterion, we check that $Y$ is nonsingular. Note that we can use  equation \ref{eq-sht-0} to solve for $b$ in terms of the two variables $\{c, d\}$, so that $\Sht^{2, \tr}(\Gamma((0)))$ is a rational surface.
\end{proof}

\subsection{Shtukas with $\mathbf{\Gamma_1(N)}$ Level Structure}\label{sec-gam-1}
Let $S = \textnormal{Spec}(A)$ be an affine scheme over $\F_q$, let $N$ be an effective divisor on $\PP^1$, and let $\tilde{\E}$ be a rank $2$ shtuka over $S$ whose pole $P$ and zero $Q$ do not belong to the support of $N \times S$.
$\tilde{\E}$ can be represented by a matrix of the form $M(T)$ as in the previous section.

A $\Gamma_1(N)$ level structure on $\tilde{\E}$ is a vector $v \in \OO_{N \times S}^2$  which generates a locally free $\OO_{N \times S}$--module of rank 1 and such that $v^\sigma = M(N)v$.
 Again, if $N$ is supported on points of degree one, we can verify this condition locally: at each point $R$ appearing in $N$ with multiplicity $m_R$, we choose the uniformizer $\pi = T-R$ at $R$ and require
\[v_R^\sigma = M(\pi + R)v_R \mod \pi^{m_R}.\]

\begin{thm}\label{thm-gam-1}
	Let $N$ be an effective divisor on $\PP^1$ of degree at least $2$. Then $\Sht^{2, \tr}(\Gamma_1(N))$ contains as a dense open substack a surface $S_N$ over $((\PP^1 \setminus N) \times (\PP^1 \setminus N)) \setminus \Delta$, where $\Delta$ is the diagonal subscheme. Moreover,
	\begin{enumerate}[\normalfont(1)]
		\item If $\deg N = 2$ and $N$ is supported on degree $1$ points, then $S_N$ is a rational surface.
	\end{enumerate}
\end{thm}
\begin{proof}
	As in the proof of Theorem \ref{thm-gam-N}, the general statement follows from the degree $2$ case, so let $N$ be a degree $2$ effective divisor on $\PP^1$. We first prove the claim in the case where $N$ is supported on degree $1$ points. Possibly after a change of variables, we can assume that $N = (0) + (\infty)$ or $N = 2(0)$. 
	
	Fix a ring $A$ over $\F_q$, $S := \text{Spec}(A)$, and two nonzero points $P, Q \in \PP^1(A)$, disjoint from $0$ and $\infty$ as necessary. An $A$--point of $\Sht^{2, \tr} (\Gamma_1(N))$ is given locally on $S$ by a pair $(M(T), v)$, where $M(T)$ is a matrix of the form
	$ M(T) := \begin{bmatrix}
	a_0 + a_1 T & b_0 + b_1 T \\
	c_0 + c_1 T & d_0 + d_1 T
	\end{bmatrix}$
	with determinant $\det M(T) = n(T-P)(T-Q)$, $n \in A^*$ , and $v$ is a nonzero vector in $\OO_{N \times S}^2$ satisfying $v^\sigma = M(N)v$. This pair is considered up to the action of $\text{GL}_2(A)$, which sends $(M(T), v)$ to $(Z^\sigma M(T)Z^{-1}, Zv)$.
	
	Suppose first that $N = (0) + (\infty)$. Restrict to the dense open substack $S_N$ of $\Sht^{2, \tr} (\Gamma_1(N))$ whose $A$--points are given locally on $S$ by pairs $(M(T), v)$ such that the restrictions of $v$ to $0$ and $\infty$ are linearly independent.
	Note that there is a unique $Z \in \text{GL}_2(A)$ which takes $v$ to the vector that is $\begin{bmatrix}
	0\\
	1
	\end{bmatrix}$ at $0$ and $\begin{bmatrix}
	1 \\ 0
	\end{bmatrix}$ at $\infty$. Hence an element of $S_N$ is represented by a unique matrix $M(T)$ as above that satisfies the equalities $M(0)\begin{bmatrix}
	0\\
	1
	\end{bmatrix} = \begin{bmatrix}
	0\\
	1
	\end{bmatrix}$ and $M(\infty)\begin{bmatrix}
	1\\
	0
	\end{bmatrix} = \begin{bmatrix}
	1\\
	0
	\end{bmatrix}$.
	That is, it can be represented by a well--defined matrix of the form
	$ M(T) := \begin{bmatrix}
	a +  T & b T \\
	c  & 1 + d T
	\end{bmatrix},$
	with $\det M(T) = a + (1 + ad - bc)T + dT^2 = n(T-P)(T-Q), n \in A^*$, which is equivalent to 
	\[
	\{
	d = n,\:
	a = nPQ,\:
	1 + ad - bc = -n(P+Q)
	\}.
	\]
	Hence $\Sht^{2, \tr}(\Gamma_1((0) + (\infty)))$ contains as a dense open substack the surface inside $\A^3_{b, c, d}$ cut out by 
	\[ d^2PQ + d(P+Q) + 1 - bc = 0. \]
	Its projective closure is nonsingular by the Jacobian criterion, and we conclude as in the proof of Theorem \ref{thm-gam-N} that it is a rational surface.
	
	Now let $N = 2(0)$, and let $\pi = T$ be a uniformizer at $0$. Restrict to the dense open substack $S_N$ of $\Sht^{2, \tr} (\Gamma_1(N))$ whose $A$--points are given locally on $S$ by pairs $(M(T), v = v_0 + v_1 \pi)$ such that $v_0$ and $v_1$ are linearly independent. There is a unique
	  $Z \in \GL_2(A)$ sending $v \in \OO_{N \times S}^2$ to $\begin{bmatrix}
	  \pi \\
	  1
	  \end{bmatrix}$; therefore, an element of $S_N(A)$ is described by a well--defined matrix $M(T)$ as above such that $M(\pi)\begin{bmatrix}
	\pi \\
	1
	\end{bmatrix} = \begin{bmatrix}
	\pi \\
	1
	\end{bmatrix} \mod \pi^2$. Comparing the coefficients of $1$ and $\pi$, this is equivalent to the system
	\[
	\{
	b_0 = 0,\:
	d_0 = 1,\:
	a_0 + b_1 = 1,\:
	c_0 + d_1 = 0
	\},
	\]
	so after renaming $a_1, b_1, c_1, d_1$ as $a,b,c,d$, we have
	$M(T) = \begin{bmatrix}
	1 - b + aT & bT \\
	-d + cT & 1 + dT
	\end{bmatrix}$. Imposing $\det M(T) = n(T^2 - (P+Q)T + PQ)$ yields the equations
	\[
	\{
	ad - bc = n,\:
	a + d = -n(P+Q),\:
	1 - b = nPQ
	\},
	\]
	and after eliminating the variables $\{a, b\}$, we obtain that this open substack of $\Sht^{2, \tr}(\Gamma_1(2(0)))$ is the surface inside $\A^3_{c, d, n}$ cut out by the equation
	\[ d^2 + dn(P+Q) + n + c - cnPQ = 0. \]
	Again, its projective closure is a nonsingular rational surface.
	
	Finally, suppose that $N = (R_0) + (R_1)$, where $R_0 \neq R_1$ are points of degree $m_0$ and $m_1$, respectively. Then for any algebraically closed field $L$ over $\F_q$, $\OO_{R_i \times \text{Spec}(L)}^2$ is isomorphic to $ (L^{2})^{m_i}$, and under this isomorphism we can identify a vector $v_i \in \OO_{R_i \times \text{Spec}(L)}^2$ with a tuple $(v_{i,1}, \dots, v_{i, m_i})$,
	
	 \noindent with $v_{i, j} \in L^2\setminus\{0\}$. Then to make sure that a generic element in $\Sht^{2, \tr}(\Gamma_1(N))$ has a well--defined representative, it suffices to observe that, generically, there is a unique matrix $Z \in \GL_2(L)$ sending $v_{0,1}$ to $\begin{bmatrix}
	0 \\ 1
	\end{bmatrix}$ and $v_{1,1}$ to $\begin{bmatrix}
	1 \\ 0
	\end{bmatrix}$. The case where $R_0 = R_1$ is analogous.
	
\end{proof}

\subsection{Shtukas with $\mathbf{\Gamma_0(N)}$ Level Structure}\label{sec-gam-0}

Let $S = \textnormal{Spec}(A)$ be an affine scheme over $\F_q$, let $N$ be an effective divisor on $\PP^1$, and let $\tilde{\E}$ be a rank $2$ shtuka over $S$ whose pole $P$ and zero $Q$ do not belong to the support of $N \times S$. 

An $A$--point of $\Sht^{2, \tr} (\Gamma_0(N))$ is represented by a pair $(M(T), v)$, where $M(T)$ is as in Section \ref{sec-gam}, and $v$ is a projective vector in $\PP^1(\OO_{N \times S})$ satisfying $v^\sigma = M(N)v$, where $M(N)$ acts on $v$ as a linear fractional transformation. This pair is considered up to the action of $\text{GL}_2(A)$, which sends $(M(T), v)$ to $(Z^\sigma M(T)Z^{-1}, Zv)$.
As in the previous section, if $N$ is supported on degree one points, the level structure condition can be checked using a uniformizer at each point in its support.

The next lemma will be used in the proof of theorem \ref{thm-gam-0}.

\begin{lemma}\label{lemm-trans}
	Let $A$ be a ring over $\F_q$.
	Fix a place $\mathit{v}$ of $\PP^1_{\F_q}$, let $\pi$ be a uniformizer at $\mathit{v}$, and denote $R_n := \F_q[\pi]/(\pi^n)$ for $n \geq 1$.  
	\begin{enumerate}[\normalfont(1)]
		\item Given any $x + y\pi \in \PP^1(A \otimes R_2)$, $z \in \PP^1(A)$ such that $x - z, y \in A^*$, there is a unique matrix $Z \in \GL_2(A)$ such that $Z(x + y\pi) = \pi \mod \pi^2$ and $Zz = \infty$.
		\item Given any $x + y\pi + z\pi^2 \in \PP^1(A \otimes R_3)$ with $y \in A^*$, there is a unique matrix $Z \in \GL_2(A)$ such that $Z(x + y\pi + z\pi^2) = \pi + \pi^2 \mod \pi^3$.
	\end{enumerate}
\end{lemma}

\begin{proof}
	\underline{Claim 1}: Fix $x + y\pi, z$ as above and let $Z := \begin{bmatrix} a & b \\ c & d
	\end{bmatrix} \in \GL_2(A)$. Then $Z(x + y\pi) = \pi \mod \pi^2$ is equivalent to 
	$ \{
	ax + b = 0,
	ay = cx + d
	\},$
	and $Zz = \infty$ is equivalent to $cz + d = 0$. Therefore both conditions can only hold simultaneously when $\dsty Z := a\begin{bmatrix} 1 & -x \\ \frac{y}{x-z} & \frac{-yz}{x-z}
	\end{bmatrix}, \text{ for } a \in A^*$.
	
	\underline{Claim 2}:  Given $x + y\pi + z\pi^2$ as above and $Z := \begin{bmatrix} a & b \\ c & d
	\end{bmatrix} \in \GL_2(A)$, the equation $Z(x + y\pi + z\pi^2) = \pi + \pi^2 \mod \pi^3$  is equivalent to 
	$ \{
	ax + b = 0,
	ay = cx + d,
	az = cy + cx + d
	\}.$ This means that $Z$ is of the form $a\begin{bmatrix} 1 & -x \\ \frac{z-y}{y} & \frac{y^2 - xz + xy}{y}
	\end{bmatrix}, \text{ for } a \in A^*$.
\end{proof}

\begin{thm}\label{thm-gam-0}
	Let $N$ be an effective divisor on $\PP^1$ of degree at least $3$. Then $\Sht^{2, \tr}(\Gamma_0(N))$ contains as a dense open substack a  stacky quotient of
	a surface $S_N$ over $((\PP^1 \setminus N) \times (\PP^1 \setminus N)) \setminus \Delta$, where $\Delta$ is the diagonal subscheme, by the trivial action of the multiplicative group $\mathbb{G}_m$. Moreover, if $N$ is supported on degree one points then
	\begin{enumerate}[\normalfont(1)]
		\item If $\deg N = 3$, then $S_N$ is a rational surface.
		\item If $\deg N = 4$, then $S_N$ is an elliptic surface.
	\end{enumerate}
\end{thm}

\begin{proof}
	As in the proofs of Theorems \ref{thm-gam-N} and \ref{thm-gam-1}, if $N = N_0 + N_1$ where $N_0$ and $N_1$ are effective divisors and we have proved the theorem for $N_0$, then the representable map $\Sht^{2, \tr}(\Gamma_0(N)) \to \Sht^{2, \tr}(\Gamma_0(N_0))$ shows that the statement is true for $N$ as well. We will prove the statements for 
	$N$ of degree $3$ or $4$ supported on degree one points; if $N$ contains a point of higher order, we can extend the argument to $N$ just like we did in the proof of Theorem \ref{thm-gam-1}.
	
	We start by proving claim (1), so suppose that $N$ is an effective divisor of degree $3$ supported on degree $1$ points. Then, after a change of variables, we may assume that $N = (0) + (1) + (\infty)$, $N= 2(0) + (1)$, or $N = 3(0)$. Throughout the proof, $A$ will denote a ring over $\F_q$, $S = \text{Spec}(A)$, and $P, Q \in \PP^1(A)\setminus  \text{Supp}(N)$ will be two distinct points.
	
	\underline{Case 1}: Let $N = (0) + (1) + (\infty)$. An $A$--point of $\Sht^{2, \tr} (\Gamma_0(N))$ is given by a pair $(M(T), v)$ as above, up to the action of $\GL_2(A)$.
	
	Note that up to multiplication by a scalar matrix, there is only one matrix $Z \in \text{GL}_2(A)$ sending $v$ to the element of $\PP^1(\OO_{N \times S})$ which is $0$ at $0$, $1$ at $1$ and $\infty$ at infinity (this is possible as long as the restrictions of $v$ to $0$, $1$ and $\infty$ are distinct, which is an open condition).
	Imposing $M(0)0 = 0$, $M(1)1 = 1$ and $M(\infty)\infty = \infty$ yields
	\[ \{
	b_0 = 0,\:
	c_1 = 0,\:
	a_0 + a_1 + b_1 = c_0 + d_0 + d_1
	\},\]
	so the element $(M(T), v)$ can be represented by a matrix
	$ M(T) := \begin{bmatrix}
	a_0 +  a_1T & b_1 T \\
	c_0  & d_0 + d_1 T
	\end{bmatrix}$
	satisfying $a_0 + a_1 + b_1 = c_0 + d_0 + d_1$ and $\det M(T) = n(T-P)(T-Q), n \in A^*$, and this representative is unique up to scalars.
	The determinant condition is equivalent to imposing
	\[
	\{
	a_0d_0 - PQa_1d_1 = 0,\:
	a_0d_1 + a_1d_0 - b_1c_0 + (P+Q)a_1d_1 = 0
	\}.
	\]
	Therefore $S_N$ is the surface inside $\PP^5_{a_0, a_1, b_1, c_0, d_0, d_1}$ described by the equations
	\[
	\{
	a_0 + a_1 + b_1 = c_0 + d_0 + d_1,\:
	a_0d_0 - PQa_1d_1 = 0,\:
	a_0d_1 + a_1d_0 - b_1c_0 + (P+Q)a_1d_1 = 0
	\}.
	\]
	
	It is not clear from this description that $S_N$ is a rational surface. To prove this claim, we start by noticing that imposing the condition that the determinant of $M(T)$ is a unit multiple of $(T-P)(T-Q)$ is equivalent to requiring that the matrix $M(T)$ has rank $1$ exactly at $T=P$ and $T=Q$. That is, there must exist a pair of points $(u, v) = ([u_0 : u_1], [v_0 : v_1])$ in $\PP^1 \times \PP^1$ such that
	\[  M(P)\begin{bmatrix}
	u_0\\u_1
	\end{bmatrix} = \begin{bmatrix}
	0\\0
	\end{bmatrix} \quad \text{and} \quad M(Q)\begin{bmatrix}
	v_0\\v_1
	\end{bmatrix} = \begin{bmatrix}
	0\\0
	\end{bmatrix},\]
	where the operation is matrix multiplication.
	
	After introducing these new variables, the equations for $S_N$ become
	\[
	\begin{cases}
	a_0 + a_1 + b_1 = c_0 + d_0 + d_1,\\
	(a_0 + a_1P)u_0 + b_1Pu_1 = 0,\\
	c_0u_0 + (d_0 + d_1P)u_1 = 0,\\
	(a_0 + a_1Q)v_0 + b_1Qv_1 = 0,\\
	c_0v_0 + (d_0 + d_1Q)v_1 = 0.
	\end{cases}
	\]
	The key fact is that we can solve for the variables $a_0, a_1, b_1, c_0, d_0, d_1$ in terms of $u_0, u_1, v_0, v_1$. Collecting coefficients we obtain the matrix
	\begin{table}[H]
		\centering
		\renewcommand{\arraystretch}{1.2}
		\begin{tabular}{cccccc}
			$a_0$ & $a_1$ & $b_1$ & $c_0$ & $d_0$ & $d_1$ \\
			\hline
			$u_0$ & $u_0P$ & $u_1P$ & $0$ & $0$ & $0$ \\
			$0$ & $0$ & $0$ & $u_0$ & $u_1$ & $u_1P$ \\
			$v_0$ & $v_0Q$ & $v_1Q$ & $0$ & $0$ & $0$ \\
			$0$ & $0$ & $0$ & $v_0$ & $v_1$ & $v_1Q$ \\
			$1$ & $1$ & $1$ & $-1$ & $-1$ & $-1$. \\
		\end{tabular}
		\renewcommand{\arraystretch}{1}
	\end{table}
	
	Taking minors of this matrix we find out that
	
	\begin{equation}\label{eq-coeffs-u-v}
	\begin{cases}
	a_0 = PQ(u_1v_0 - u_0v_1)((P-Q)u_1v_1 + (Q-1)u_0v_1 - (P-1)u_1v_0),\\
	a_1 = (u_1v_0P - u_0v_1Q)((P-Q)u_1v_1 + (Q-1)u_0v_1 - (P-1)u_1v_0),\\
	b_1 =  (P-Q)u_0v_0((P-Q)u_1v_1 + (Q-1)u_0v_1 - (P-1)u_1v_0),\\
	c_0 = (P-Q)u_1v_1((Q-P)u_0v_0 + P(1 - Q)u_1v_0 + Q(P-1)u_0v_1),\\
	d_0 = (u_1v_0P - u_0v_1Q)((Q-P)u_0v_0 + P(1 - Q)u_1v_0 + Q(P-1)u_0v_1),\\
	d_1 = (u_1v_0 - u_0v_1)((Q-P)u_0v_0 + P(1 - Q)u_1v_0 + Q(P-1)u_0v_1).
	\end{cases}
	\end{equation}
	
	This shows that $S_N$ is a rational surface. Note that each of the $a_0, \dots, d_1$ is a $(2, 2)$--form in the variables $(u, v)$.

	\underline{Case 2}: Let $N = 2(0) + (\infty)$, and let $(M(T), v)$ represent an $A$--point of $\Sht^{2, \tr} (\Gamma_0(N))$. Let $\pi = T$ be a uniformizer at $0$.
	
	By Lemma \ref{lemm-trans}, generically there is a unique matrix $Z \in \text{GL}_2(A)$ up to scalar multiplication sending $v$ to the element of $\PP^1(\OO_{N \times S})$ which is $\pi$ at $2(0)$ and $\infty$ at infinity.
	Imposing the conditions $M(\pi)\pi = \pi \mod \pi^2$, $M(\infty) = \infty$ and and $\det M(T) = n(T-P)(T-Q)$ for some $n \in A^*$, we find that 
	$S_N$ is the surface inside $\PP^5_{a_0, a_1, b_1, c_0, d_0, d_1}$ described by the equations
	\[
	\{
	a_0 + b_1 = d_0,\:
	a_0d_0 - PQa_1d_1 = 0,\:
	a_0d_1 + a_1d_0 - b_1c_0 + (P+Q)a_1d_1 = 0
	\}.
	\]
	As in case 1, we can solve for $a_0, \dots, d_1$ in terms of the variables $(u, v)$, so $S_N$ is again a rational surface.
	
	\underline{Case 3}: Let $N = 3(0)$, let $(M(T), v)$ represent an $A$--point of $\Sht^{2, \tr} (\Gamma_0(N))$, and let $\pi = T$ be a uniformizer at $0$.
	
	Lemma \ref{lemm-trans} shows that, generically, there is a unique matrix $Z \in \text{GL}_2(A)$ up to scalars sending $v$ to the element of $\PP^1(\OO_{N \times S})$ which is $\pi + \pi^2$ at $3(0)$.
	Imposing the conditions $M(\pi)(\pi + \pi^2) = \pi + \pi^2 \mod \pi^3$ and $\det M(T) = n(T-P)(T-Q)$ for some $n \in A^*$, we find that 
	$S_N$ is the surface inside $\PP^6_{a_0, a_1, b_1, c_0, c_1, d_0, d_1}$ described by the equations
	\begin{align*}
	\{
	&a_0 + b_1 = d_0,\:
	a_0 + a_1 = c_0 + d_0 + d_1,\:
	a_0d_0 - PQ(a_1d_1 - b_1c_1) = 0,\\
	&a_0d_1 + a_1d_0 - b_1c_0 + (P+Q)(a_1d_1 - b_1c_1) = 0
	\}.
	\end{align*}
	We conclude as in the previous cases that $S_N$ is a rational surface.
	
	Now we will prove claim (2). After performing an appropriate change of variables, there are five cases to consider: $N = (0) + (1) + (\infty) + (R)$ for $R \in \PP^1\setminus\{0, 1, \infty\}$, $N = 2(0) + (1) + (\infty)$, $N = 2(0) + 2(\infty)$, $N = 3(0) + (\infty)$, and $N = 4(0)$.
	
	\underline{Case 1}: Let $N = (0) + (1) + (\infty) + (R)$. To get the equations for $S_N$, we start from those of $S_{N_0}$ with $N_0 = (0) + (1) + (\infty)$ and add the condition that $M(R)w = w^q$ for some projective vector $w$. That is, $S_N$ is the hypersurface of $\PP^1_{u} \times \PP^1_v \times \PP^1_{w}$ given by
	\[ c_0w^{q+1} + (d_0 + d_1R)w^q - (a_0 + a_1R)w - b_1R = 0, \]
	where the variables $a_0, \dots, d_1$ are defined in terms of $u$ and $v$ as in (\ref{eq-coeffs-u-v}). Therefore $S_N$ is a $(2, 2, q+1)$--surface inside $\PP^1_{u} \times \PP^1_v \times \PP^1_{w}$. Recall that a smooth $(2,2)$--curve in $\PP^1_{u} \times \PP^1_v$ is an elliptic curve; hence the projection onto the $w$ line shows that $S_N$ is an elliptic surface.
	
	\underline{Case 2}: Let $N = 2(0) + (1) + (\infty)$. We start from the equations for $S_{(0) + (1) + (\infty)}$ and add the condition $M(\pi)w\pi = w^q\pi \mod \pi^2$, where $\pi$ is a uniformizer at $0$; this shows that $S_N$ is the surface inside $\PP^5_{a_0, a_1, b_1, c_0, d_0, d_1} \times \PP^1$ described by the equations
	\[
	\{
	d_0w^q - a_0w - b_1 = 0,\:
	a_0 + a_1 + b_1 = c_0 + d_0 + d_1,\:
	a_0d_0 - PQa_1d_1 = 0,\:
	a_0d_1 + a_1d_0 - b_1c_0 + (P+Q)a_1d_1 = 0
	\}.
	\]
	Again, we can write everything in terms of the variables $u, v, w$ to show that $S_N$ is an elliptic surface.
	
	\underline{Case 3}: Let $N = 2(0) + 2(\infty)$. We impose the conditions $M(\pi_0)\pi_0 = \pi_0 \mod \pi_0^2$ and $M(\pi_\infty + R)(1 + w\pi_\infty) =  1 + w^q\pi_\infty \mod \pi_\infty^2$, where $\pi_0 = T$ and $\pi_\infty = \frac{1}{T}$ are uniformizers at $0$ and $\infty$, respectively. Then $S_N$ is the surface inside $\PP^6_{a_0, a_1, b_1, c_0, c_1 d_0, d_1} \times \PP^1_w$ described by
	\begin{align*}
	\{
	&(c_1 + d_1)w^q + (c_1 - a_1)w + c_0 + d_0 - a_0 = 0,\:
	a_0 + b_1 = d_0,\:
	a_1 + b_1 = c_1 + d_1,\\
	&a_0d_0 - PQ(a_1d_1 - b_1c_1) = 0,\:
	a_0d_1 + a_1d_0 - b_1c_0 + (P+Q)(a_1d_1 - b_1c_1) = 0
	\}.
	\end{align*}
	We rewrite $a_0, \dots, d_1$ in terms of the variables $u, v, w$ to show that $S_N$ is an elliptic surface.
	
	\underline{Case 4}: Let $N = 3(0) + (\infty)$. Then after adding the new condition at $\infty$, we get that $S_N$ is the elliptic surface inside $\PP^6_{a_0, a_1, b_1, c_0, c_1, d_0, d_1} \times \PP^1_w$ described by the equations
	\begin{align*}
	\{
	&c_1w^{q+1} + d_1w^q - a_1w - b_1 = 0,\:
	a_0 + b_1 = d_0,\:
	a_0 + a_1 = c_0 + d_0 + d_1,\\
	&a_0d_0 - PQ(a_1d_1 - b_1c_1) = 0,\:
	a_0d_1 + a_1d_0 - b_1c_0 + (P+Q)(a_1d_1 - b_1c_1) = 0
	\}.
	\end{align*}
	
	\underline{Case 5}: Finally, let $N = 4(0)$. Then to describe $S_N$ we must impose the condition $M(\pi_0)(\pi_0 + \pi_0^2 + w\pi_0^3) = (\pi_0 + \pi_0^2 + w\pi_0^3) \mod \pi_0^4$. This yields that
	$S_N$ is the elliptic surface inside $\PP^6_{a_0, a_1, b_1, c_0, c_1, d_0, d_1} \times \PP^1_w$ described by the equations
	\begin{align*}
	\{
	&d_0 w^q - a_0 w -a_1 + 2c_0 + c_1 + d_1 = 0,\:
	a_0 + b_1 = d_0,\:
	a_0 + a_1 = c_0 + d_0 + d_1,\\
	&a_0d_0 - PQ(a_1d_1 - b_1c_1) = 0,\:
	a_0d_1 + a_1d_0 - b_1c_0 + (P+Q)(a_1d_1 - b_1c_1) = 0
	\}.
	\end{align*}
\end{proof}

\section{$\mathbf{\textbf{Sht}^{2, \textbf{tr}}(\Gamma_0(N))}$ for $\mathbf{N}$ of Degree Four}\label{sec-deg-4}
\graphicspath{{1_Body/Figures/}}
In this section we study $\Sht^{2,\tr}(\Gamma_0(N))$  for 
$N$ an effective divisor on $X := \PP^1_{\F_q}$ of degree $4$ such that every point in the support of $N$ has degree $1$. In this case, $\Sht^{2,\tr}(\Gamma_0(N))$ is an elliptic surface whose compactification $\overline{\Sht}^{2, \tr}(\Gamma_0(N))$ can be described as a hypersurface in $\PP^1 \times \PP^1 \times \PP^1$, possibly with singularities.
Hence Lemma \ref{lemm-arith-gen} gives an upper bound on the arithmetic genus of $\overline{\Sht}^{2, \tr}(\Gamma_0(N))$.

\subsection{$\mathbf{N = 4(0)}$}\label{sec-4-0}
Let $N$ be a degree $4$ divisor on $\PP^1_{\F_q}$ supported on a degree $1$ point; after a change of variables, we may assume that $N = 4(0)$. Then the proof of Theorem \ref{thm-gam-0} shows that we can find a description of $\overline{\Sht}^{2, \tr}(\Gamma_0(4(0)))$ as a hypersurface of $\PP^1 \times \PP^1 \times \PP^1$ given by a polynomial of multidegree $(2, 2, q)$. This surface is singular, with a one dimensional singular subscheme. By Lemma \ref{lemm-arith-gen}, we can conclude that the arithmetic genus of $\overline{\Sht}^{2, \tr}(\Gamma_0(4(0)))$ is less than or equal to $q$. However, experimental data suggests that this upper limit is never reached.

We used Magma to compute several invariants of $\overline{\Sht}^{2, \tr}(\Gamma_0(4(0)))$ for small values of $q$, which we present in Table \ref{tab-4-0}. Here $b_2$ denotes the second Betti number of the surface, rk(T) denotes the rank of its trivial lattice, $\#$BF the number of singular fibers, and ``Types'' the Kodaira types of the singular fibers, using the notation $(\text{K}, n)$ to represent a place of bad reduction of degree $n$ whose corresponding fiber is of Kodaira type K. The last two rows collect our conjectured values for these invariants, based on the presented data.

\begin{table}[H]
	\centering
	\renewcommand{\arraystretch}{1.2}
	\begin{tabular}{|c|c|c|c|c|}
		\hline
		$\mathbf{q}$ &  $\mathbf{b_2}$ & \textbf{rk(T)} &$\#$\textbf{BF} & \textbf{Types}\\
		\hline
		2 & 10 & 8 & 2 & $(\text{III}, 1), (\text{I}_6, 1)$ \\
		\hline
		4 & 22 & 18 & 2 & $(\text{III}^*, 1), (\text{I}_{10}, 1)$ \\
		\hline
		8 & 34 & 26 & 2 &  $(\text{III}^*, 1), (\text{I}_{18}, 1)$ \\
		\hline
		16 & 58 & 42 & 2 &  $(\text{III}^*, 1), (\text{I}_{34}, 1)$ \\
		\hline
		3 & 22 & 17  & 8 & $(\text{I}_4^*, 1), (\text{I}_8, 1), (\text{I}_1, 3), (\text{I}_1, 3)$ \\
		\hline
		5 & 34 & 25 & 12 & $(\text{I}_8^*, 1), (\text{I}_{12}, 1), (\text{I}_1, 5), (\text{I}_1, 5)$  \\
		\hline
		7 & 46 & 33 & 16 & $(\text{I}_{12}^*, 1), (\text{I}_{16}, 1), (\text{I}_1, 7), (\text{I}_1, 7)$  \\
		\hline
		9 & 58 & 41 & 20 & $(\text{I}_{16}^* 1), (\text{I}_{20}, 1), (\text{I}_1, 9), (\text{I}_1, 9)$ \\
		\hline
		11 & 70 & 49 & 24 & $(\text{I}_{20}^*, 1), (\text{I}_{24}, 1), (\text{I}_1, 11), (\text{I}_1, 11)$ \\
		\hline
		even $q > 2$ & $3q + 10$ & $2q + 10$ & 2 &  $(\text{III}^*, 1), (\text{I}_{2q+2}, 1)$ \\
		\hline
		odd $q$ &  $6q + 4$ &  $4q + 5$ & $2q + 2$ & $(\text{I}_{2q-2}^*, 1), (\text{I}_{2q+2}, 1), (\text{I}_1, q), (\text{I}_1, q)$  \\
		\hline
	\end{tabular}
	\renewcommand{\arraystretch}{1}
	\caption{Invariants of $\overline{\Sht}^{2, \tr}(\Gamma_0(4(0)))$}
	\label{tab-4-0}
\end{table}

In particular, by checking the Betti numbers we conclude that $\overline{\Sht}^{2, \tr}(\Gamma_0(4(0)))$ is a rational elliptic surface when $q=2$ and a K3 elliptic surface when $q=3$ or $q=4$.

\begin{conj}
	Let $q$ be a prime power and let $S_q := \overline{\Sht}^{2, \tr}(\Gamma_0(4(0)))$ over $F := \F_q(P,Q)$.
	
	\noindent If $q$ is even, then
	\begin{enumerate}[\normalfont(1)]
		\item $p_a(S_q) = \left[ \frac{q}{4} \right] + 1$ 
		and $b_2(S_q) = 12\left[ \frac{q}{4} \right] + 10$.
		\item The trivial lattice of $S_q$ has rank $8$ when $q = 2$ and rank $2q + 10$ for $q > 2$.
		\item $S_q$ has two singular fibers, one of type $\textnormal{I}_{2q + 2}$ at $1$, and one of type $\textnormal{III}^*$ at $\infty$ (except when $q=2$, for which the fiber at $\infty$ is of type $\textnormal{III}$).
	\end{enumerate}
	If $q$ is odd, then
	\begin{enumerate}[\normalfont(1)]
		\item $p_a(S_q) = \frac{q+1}{2}$ 
		and $b_2(S_q) = 6q + 4$.
		\item The trivial lattice of $S_q$ has rank $4q + 5$.
		\item $S_q$ has $2q + 2$ singular fibers over $\overbar{F}$: one of type $\textnormal{I}_{2q-2}^*$ at $\infty$, one of type $\textnormal{I}_{2q + 2}$ at $1$, and the rest coming from two bad fibers of type $\textnormal{I}_1$, of degree $q$ over $F$.
	\end{enumerate}
\end{conj}

\subsection{$\mathbf{N = 2(0) + 2(\infty)}$}\label{sec-2-2}
Let $N = 2(A) + 2(B)$ be a degree $4$ divisor on $\PP^1_{\F_q}$, where $A$ and $B$ are degree 1 points. Without loss of generality, we may assume $N = 2(0) + 2(\infty)$. We know that 
$\overline{\Sht}^{2, \tr}(\Gamma_0(2(0) + 2(\infty)))$ is birational to a hypersurface of $\PP^1 \times \PP^1 \times \PP^1$  of multidegree $(2, 2, q)$. This surface has a one dimensional singular subscheme, except when $q=2$, when it is smooth. Hence
by Lemma \ref{lemm-arith-gen} its arithmetic genus is at most $q$, and it is exactly $2$ for $q=2$. The experimental data collected in Table \ref{tab-2-2} suggests that the actual arithmetic genus is generally smaller.

\begin{table}[H]
	\centering
	\renewcommand{\arraystretch}{1.2}
	\begin{tabular}{|c|c|c|c|c|}
		\hline
		$\mathbf{q}$ &  $\mathbf{b_2}$ & \textbf{rk(T)} &$\#$\textbf{BF} & \textbf{Types}\\
		\hline
		2 & 22 & 18 & 3 & $(\text{I}_2^*, 1), (\text{I}_6, 1), (\text{I}_6, 1)$ \\
		\hline
		4 & 34 & 26 & 3 & $(\text{I}_2^*, 1), (\text{I}_{10}, 1), (\text{I}_{10}, 1)$  \\
		\hline
		8 & 58 & 42 & 3 &  $(\text{I}_2^*, 1), (\text{I}_{18}, 1), (\text{I}_{18}, 1)$ \\
		\hline
		3 & 22 & 17  & 9 & $(\text{I}_2, 1), (\text{I}_8, 1), (\text{I}_8, 1), (\text{I}_1, 2), (\text{I}_1, 4)$ \\
		\hline
		5 & 34 & 25 & 13 & $(\text{I}_2, 1), (\text{I}_{12}, 1), (\text{I}_{12}, 1), (\text{I}_1, 4), (\text{I}_1, 6)$  \\
		\hline
		7 & 46 & 33 & 17 & $(\text{I}_2, 1), (\text{I}_{16}, 1), (\text{I}_{16}, 1), (\text{I}_1, 6), (\text{I}_1, 8)$  \\
		\hline
		9 & 58 & 41 & 21 & $(\text{I}_2, 1), (\text{I}_{20}, 1), (\text{I}_{20}, 1), (\text{I}_1, 8), (\text{I}_1, 10)$ \\
		\hline
		11 & 70 & 49 & 25 & $(\text{I}_2, 1), (\text{I}_{24}, 1), (\text{I}_{24}, 1) (\text{I}_1, 10), (\text{I}_1, 12)$ \\
		\hline
		even $q$ & $6q + 10$ &  $4q  + 10$ & 3 & $(\text{I}_2^*, 1), (\text{I}_{2q + 2}, 1), (\text{I}_{2q + 2}, 1)$ \\
		\hline
		odd $q$  & $6q + 4$ & $4q + 5$ & $2q + 3$ & $(\text{I}_2, 1), (\text{I}_{2q + 2}, 1), (\text{I}_{2q + 2}, 1), $ \\
		
		&  &  &  & $ (\text{I}_1, q - 1), (\text{I}_1, q + 1)$ \\
		\hline
	\end{tabular}
	\renewcommand{\arraystretch}{1}
	\caption{Invariants of $\overline{\Sht}^{2, \tr}(\Gamma_0(2(0)+2(\infty)))$}
	\label{tab-2-2}
\end{table}

Note that $\overline{\Sht}^{2, \tr}(\Gamma_0(2(0)+2(\infty))$ is a K3 surface when $q=2$ or $q=3$.

\begin{conj}
	Let $q$ be a prime power and let $S_q$ be $\overline{\Sht}^{2, \tr}(\Gamma_0(2(0) + 2(\infty)))$ over $F := \F_q(P, Q)$.
	
	\noindent If $q$ is even, then
	\begin{enumerate}[\normalfont(1)]
		\item $p_a(S_q) =  \frac{q+2}{2}$ 
		and $b_2(S_q) = 6q + 10$.
		\item The trivial lattice of $S_q$ has rank $4q + 10$.
		\item $S_q$ has two singular fibers of type $\textnormal{I}_{2q + 2}$ at $0$ and $1$, and one of type $\textnormal{I}_2^*$ at $\infty$.
	\end{enumerate}
	If $q$ is odd, then
	\begin{enumerate}[\normalfont(1)]
		\item $p_a(S_q) = \frac{q+1}{2}$ 
		and $b_2(S_q) = 6q + 4$.
		\item The trivial lattice of $S_q$ has rank $4q + 5$.
		\item $S_q$ has $2q + 3$ singular fibers over $\overbar{F}$: two of type $\textnormal{I}_{2q + 2}$ at $0$ and $-1$, one of type $\textnormal{I}_2$ at $\infty$, and the rest corresponding to two bad fibers of type $\textnormal{I}_1$ over $F$, of degrees $q-1$ and $q+1$, respectively.
	\end{enumerate}
\end{conj}

\subsection{$\mathbf{N = 3(0)+(\infty)}$}\label{sec-3-1}
Let $N = 3(A) + (B)$ be a divisor on $\PP^1_{\F_q}$, where $A$ and $B$ are degree 1 points; we may assume $N = 3(0) + (\infty)$. Up to birational equivalence, we can describe 
$\overline{\Sht}^{2, \tr}(\Gamma_0(3(0) + (\infty)))$ as a hypersurface cut out by a polynomial of multi--degree $(2, 2, q + 1)$ inside $\PP^1 \times \PP^1 \times \PP^1$, whose singular subscheme has dimension zero. Therefore the arithmetic genus of $\overline{\Sht}^{2, \tr}(\Gamma_0(3(0) + (\infty)))$ is at most $q+1$, and based on the  experimental data, we expect it to equal $q$.

We need to make a remark about the way we collected our data. When an elliptic surface is defined over a function field $k(X)$ where $k$ is not perfect and has characteristic 2 or 3, Szydlo showed that new types of Kodaira fibers can appear \cite{Szy}; these new fibers are not supported on Magma, which uses the classical Tate algorithm. Therefore when we try to compute local information of the generic fiber of
$\overline{\Sht}^{2, \tr}(\Gamma_0(N))$, we may get an error due to the fact that we are working over the imperfect residue field $\F_q(P, Q)$, if $\overline{\Sht}^{2, \tr}(\Gamma_0(N))$ has a singular fiber which is not in Kodaira's original list. Another reason why it might not be possible to work over $\F_q(P, Q)(T)$ is that some of the computations are too slow.

While these issues did not affect our computations in Sections \ref{sec-4-0} and \ref{sec-2-2}, they do now; we are only able to compute data about the generic fiber when we work in characteristic greater than $3$. In characteristics 2 and 3, we instead choose concrete values for the pole and zero of the shtuka and compute the local information associated to the corresponding closed fiber of $\overline{\Sht}^{2, \tr}(\Gamma_0(3(0) + (\infty)))$. This is the same strategy that we will follow in Sections \ref{sec-2-1-1} and \ref{sec-mult-1}, for all characteristics.

The data we collected on $\overline{\Sht}^{2, \tr}(\Gamma_0(3(0) + (\infty)))$ is split between Tables \ref{tab-3-1-p1} and \ref{tab-3-1-p2}. 
Table \ref{tab-3-1-p2} is formatted as Table \ref{tab-4-0}, but on Table \ref{tab-3-1-p1} there are two differences. First, we introduce two new columns $P$ and $Q$ for the specific values of pole and zero that were used in the computations; there $\alpha_r$ denotes a generator of $\F_r^*$. Note that while $b_2$ stays the same if the values of $P$ and $Q$ are modified, this needs not be true for the other invariants. Secondly, in the ``Types'' column the notation K$\times n$ represents $n$ fibers of Kodaira type K over the algebraic closure, and we omit $n$ when it equals one. The tables in the next two sections will be formatted as Table \ref{tab-3-1-p1}.

\begin{table}[H]
	\centering
	\renewcommand{\arraystretch}{1.2}
	\begin{tabular}{|c|c|c|c|c|c|c|}
		\hline
		$\mathbf{q}$ &\textbf{P} & \textbf{Q} &  $\mathbf{b_2}$ & \textbf{rk(T)} &$\#$\textbf{BF} & \textbf{Types}\\
		\hline
		2 & 1 & $\alpha_4$ & 22 & 16 & 4 & $ \text{II} \times 2, \text{I}_6, \text{I}_{10} $ \\
		\hline
		4 & 1 & $\alpha_4$ & 46 & 32  & 6 &  $ \text{II} \times 4, \text{I}_{10}, \text{I}_{22} $ \\
		\hline
		8 & 1 & $\alpha_8$ & 94 & 64 & 10 &  $ \text{II} \times 8, \text{I}_{18}, \text{I}_{46} $ \\
		\hline
		3 & 1 & 2 & 34 & 24  &  14 & $\text{I}_8, \text{I}_{16}, \text{I}_1 \times 12$ \\
		\hline
		9 & 1 & 2 & 106 & 72 & 38 &  $\text{I}_{20}, \text{I}_{52}, \text{I}_1 \times 36 $  \\
		\hline
		even $q$ &  &  & $12q - 2$ & $8q$ & $q + 2$ &   $ \text{II} \times q, \text{I}_{2q+2}, \text{I}_{6q-2} $ \\
		\hline
		odd $q$ &  &  & $12q - 2$ & $8q$ & $4q + 2$ & $\text{I}_{2q + 2}, \text{I}_{6q - 2},  \text{I}_1 \times 4q $  \\
		\hline
	\end{tabular}
	\renewcommand{\arraystretch}{1}
	\caption{Invariants of $\overline{\Sht}^{2, \tr}(\Gamma_0(3(0) + (\infty)))$, part 1}
	\label{tab-3-1-p1}
\end{table}

\begin{table}[H]
	\centering
	\renewcommand{\arraystretch}{1.2}
	\begin{tabular}{|c|c|c|c|c|}
		\hline
		$\mathbf{q}$ &  $\mathbf{b_2}$ & \textbf{rk(T)} &$\#$\textbf{BF} & \textbf{Types}\\
		\hline
		5 & 58 & 40 & 22 & $(\text{I}_{12}, 1), (\text{I}_{28}, 1), (\text{I}_1, 10), (\text{I}_1, 10)$  \\
		\hline
		7 & 82 & 56 & 30 & $(\text{I}_{16}, 1), (\text{I}_{40}, 1), (\text{I}_1, 14), (\text{I}_1, 14)$  \\
		\hline
		11 & 130 & 88 & 46 & $ (\text{I}_{24}, 1), (\text{I}_{64}, 1) (\text{I}_1, 22), (\text{I}_1, 22)$ \\
		\hline
		$q$  & $12q - 2$ & $8q$ & $4q + 2$ & $(\text{I}_{2q + 2}, 1), (\text{I}_{6q - 2}, 1) (\text{I}_1, 2q), (\text{I}_1, 2q)$\\
		\hline
	\end{tabular}
	\renewcommand{\arraystretch}{1}
	\caption{Invariants of $\overline{\Sht}^{2, \tr}(\Gamma_0(3(0) + (\infty)))$, part 2}
	\label{tab-3-1-p2}
\end{table}

In particular $\overline{\Sht}^{2, \tr}(\Gamma_0(3(0) + (\infty)))$ is a K3 surface when $q=2$.

\begin{conj}
	Let $q$ be a prime power and let $S_q := \overline{\Sht}^{2, \tr}(\Gamma_0(3(0) + (\infty)))$ over $F := \F_q(P,Q)$. Then
	\begin{enumerate}[\normalfont(1)]
		\item $p_a(S_q) = q$ 
		and $b_2(S_q) = 12q-2$.
		\item The trivial lattice of $S_q$ has rank $8q$.
		\item If $q$ is even, then $S_q$ has $q + 2$ singular fibers over $\overbar{F}$: one of type $\textnormal{I}_{2q + 2}$ at $1$, one of type $\textnormal{I}_{6q - 2}$ at $0$, and $q$ of type $\textnormal{II}$.
		\item If $q$ is odd, then $S_q$ has $4q + 2$ singular fibers over $\overbar{F}$: one of type $\textnormal{I}_{2q + 2}$ at $-1$, one of type $\textnormal{I}_{6q - 2}$ at $0$, and $4q$ of type $\textnormal{I}_1$.
		
	\end{enumerate}
\end{conj}

\subsection{$\mathbf{N = 2(0)+(1)+(\infty)}$}\label{sec-2-1-1}
\begin{prop}
	Let $q$ be a prime power and let $S_q := \overline{\Sht}^{2, \tr}(\Gamma_0(2(0) + (1) + (\infty)))$ over $F := \F_q(P,Q)$. Then
	$p_a(S_q) = q$ and $b_2(S_q) = 12q-2$.
\end{prop}
\begin{proof}
	The proof of Theorem \ref{thm-gam-0} shows that 
	$S_q$ is birational to a hypersurface of the product $\PP^1_{u_0, u_1} \times \PP^1_{v_0, v_1} \times \PP^1_{w_0, w_1}$ described by a polynomial of multidegree $(2, 2, q)$. This surface has exactly eight rational double points:
	\begin{align*}
	\{&(0,1,0,1,0,1),  (P,1,P,1,1,1), (Q,1,Q,1,1,1), (P,1,Q,1,1,1), (0,1,0,1,1,0),\\ &(P,1,Q,1,1,0), (Q - P,Q-1,0,1,0,1), (0,1,P-Q,P-1,0,1)\},
	\end{align*}
	and no other singularities.
	Hence Lemma \ref{lemm-arith-gen} implies that $p_a(S_q) = q$ and $b_2(S_q) = 12q-2$.
\end{proof}

We record in Table \ref{tab-2-1-1} experimental data which allows us to predict the structure of the singular fibers of $S_q$.

\begin{table}[H]
	\centering
	\renewcommand{\arraystretch}{1.2}
	\begin{tabular}{|c|c|c|c|c|c|c|}
		\hline
		$\mathbf{q}$ &\textbf{P} & \textbf{Q} &  $\mathbf{b_2}$ & \textbf{rk(T)} &$\#$\textbf{BF} & \textbf{Types}\\
		\hline
		2 & $\alpha_8$ & $\alpha_8^3$ & 22 & 15  & 5 & $\text{II} \times 2, \text{I}_{6} \times 2, \text{I}_{4} $  \\
		\hline
		4 & $\alpha_{16}$ & $\alpha_{16}^2$ & 46 & 31  & 7 & $\text{II} \times 4, \text{I}_{10} \times 2, \text{I}_{12} $ \\
		\hline
		8 & $\alpha_8$ & $\alpha_8^3$ & 94 & 63 & 11 &  $\text{II} \times 8, \text{I}_{18} \times 2, \text{I}_{28} $ \\
		\hline
		3 & 2 & $\alpha_9$ & 34 & 23 & 15 & $\text{I}_{8} \times 2, \text{I}_{8}, \text{I}_1 \times 12 $\\
		\hline
		5 & 2 & 3 & 58 & 39 & 23 & $\text{I}_{12} \times 2, \text{I}_{16},  \text{I}_1 \times 20 $ \\
		\hline
		7 & 2 & 3 & 82 & 55  & 31 & $\text{I}_{16} \times 2, \text{I}_{24},  \text{I}_1 \times 28 $ \\
		\hline
		9 & 2 & $\alpha_9$ & 106 & 71   & 39 & $\text{I}_{20} \times 2, \text{I}_{32},  \text{I}_1 \times 36 $ \\
		\hline
		11 & 2 & 3 & 130 & 87 &  47 & $\text{I}_{24} \times 2, \text{I}_{40},  \text{I}_1 \times 44 $  \\
		\hline
		even $q$ &  &  & $12q - 2$ & $8q - 1$& $q + 3$& $\text{II} \times q, \text{I}_{2q + 2} \times 2, \text{I}_{4q - 4} $ \\
		\hline
		odd $q$ &  &  & $12q - 2$ & $8q - 1$& $4q + 3$ & $\text{I}_{2q + 2} \times 2, \text{I}_{4q - 4},  \text{I}_1 \times 4q $ \\
		\hline
	\end{tabular}
	\renewcommand{\arraystretch}{1}
	\caption{Invariants of $\overline{\Sht}^{2, \tr}(\Gamma_0(2(0) + (1) + (\infty)))$}
	\label{tab-2-1-1}
\end{table}
Notice that $\overline{\Sht}^{2, \tr}(\Gamma_0(2(0) + (1) (\infty)))$ is a K3 surface when $q=2$.

\begin{conj}
	Let $q$ be a prime power and let $S_q := \overline{\Sht}^{2, \tr}(\Gamma_0(2(0) + (1) + (\infty)))$ over $F := \F_q(P,Q)$. Then
	\begin{enumerate}[\normalfont(1)]
		
		\item The trivial lattice of $S_q$ has rank $8q - 1$.
		\item If $q$ is even, then $S_q$ has $q + 3$ singular fibers over $\overbar{F}$: two of type $\textnormal{I}_{2q + 2}$ at $0$ and $1$, one of type $\textnormal{I}_{4q - 4}$ at $\infty$, and $q$ of type $\textnormal{II}$.
		\item If $q$ is odd, then $S_q$ has $4q + 3$ singular fibers over $\overbar{F}$: two of type $\textnormal{I}_{2q + 2}$ at $0$ and $1$, one of type $\textnormal{I}_{4q - 4}$ at $\infty$, and $4q$ of type $\textnormal{I}_1$.
		
	\end{enumerate}
\end{conj}

\subsection{Multiplicity one}\label{sec-mult-1}
In this section we treat the remaining case, in which $N$ is a divisor on $\PP^1$ consisting of four distinct degree 1 points appearing with multiplicity one. If $q$ is odd we take $N = (0) + (1) + (-1) + (\infty)$, and if $q > 2$ is even we take $N = (0) + (1) + (\alpha_q) + (\infty)$, where $\alpha_q$ generates $\F_q^*$. We found a birational model for $\overline{\Sht}^{2, \tr}(\Gamma_0(N))$ which is a hypersurface of $\PP^1 \times \PP^1 \times \PP^1$ cut down by a polynomial of multi--degree $(2, 2, q+1)$; the only singularities on this model are rational double points.
We present the invariants of $\overline{\Sht}^{2, \tr}(\Gamma_0(N))$ for small values of $q$ in Table \ref{tab-mult-1}.
\begin{table}[H]
	\centering
	\renewcommand{\arraystretch}{1.2}
	\begin{tabular}{|c|c|c|c|c|c|c|}
		\hline
		$\mathbf{q}$ &\textbf{P} & \textbf{Q} &  $\mathbf{b_2}$ & \textbf{rk(T)} &$\#$\textbf{BF} & \textbf{Types}\\
		\hline
		4 & $\alpha_{16}$ & $\alpha_{16}^2$ & 58 & 38 &  9 & $\text{II} \times 5, \text{I}_{10} \times 4$\\
		\hline
		8 & $\alpha_{8}^2$ & $\alpha_{8}^3$ & 106 & 70 & 14 & $\text{III}, \text{II} \times 8, \text{I}_2, \text{I}_{16},  \text{I}_{18} \times 3 $  \\
		\hline
		3 & $\alpha_9$ & $\alpha_9^2$ & 46 & 30 &  20 & $\text{I}_{8} \times 4, \text{I}_1 \times 16$  \\
		\hline
		5 & 2 & $\alpha_{25}$ & 70 & 46 & 28 &  $\text{I}_{12} \times 4, \text{I}_1 \times 24$   \\
		\hline
		7 & 2 & 4 & 94 & 62 & 36 & $\text{I}_{16} \times 4, \text{I}_1 \times 32$ \\
		\hline
		9 & $\alpha_9$ & $\alpha_9^2$ & 118 & 78 & 44 & $\text{I}_{20} \times 4, \text{I}_1 \times 40$    \\
		\hline
		11 & 2 & 4 & 142 & 94 & 52 & $\text{I}_{24} \times 4, \text{I}_1 \times 48$    \\
		\hline
		even $q>2$ &  &  & $12q + 10$ & $8q + 6$ & &   \\
		\hline
		odd $q$ &  &  & $12q + 10$ & $8q + 6$ & $4q + 8$ & $\text{I}_{2q+2} \times 4, \text{I}_1 \times (4q + 4)$  \\
		\hline
	\end{tabular}
	\renewcommand{\arraystretch}{1}
	\caption{Invariants of $\overline{\Sht}^{2, \tr}(\Gamma_0(N))$, for $N$ a squarefree divisor}
	\label{tab-mult-1}
\end{table}
Note that $\overline{\Sht}^{2, \tr}(\Gamma_0(N))$ is never a rational or K3 surface in this case.

\begin{prop}
	Let $q$ be a prime power and let $S_q := \overline{\Sht}^{2, \tr}(\Gamma_0(N))$ over the field $F := \F_q(P,Q)$, where $N = (0) + (1) + (-1) + (\infty)$ if $q$ is odd and $N = (0) + (1) + (\alpha_q) + (\infty)$ if $q>2$ is even. Then
	 $p_a(S_q) = q + 1$ and $b_2(S_q) = 12q +10$.
\end{prop}
\begin{proof}
From the proof of Theorem \ref{thm-gam-0}, we have that $S_q$ is birational to a hypersurface of $\PP^1 \times \PP^1 \times \PP^1$ cut down by a polynomial of multi--degree $(2, 2, q+1)$. If $q$ is odd, the only singularities of $S_q$ are
\begin{align*}
	\{ &(0,1,0,1,0,1), (0,1,0,1,-1,1), (1,0,1,0,1,0), (1,0,1,0,-1,1),
	 (0,1,P - Q,P-1,0,1), \\&(Q -P, Q-1,0,1,0,1), (P-PQ,P-Q,1,0,1,0), (1,0,Q-QP,Q-P,1,0) \}
\end{align*}
and for even $q \geq 4$ the singularities are
\begin{align*}
\{ &(0,1,0,1,0,1), (1,1,1,1,1,1), (P,Q,1,1,1,1),(1,1,Q,P,1,1),
(Q,\alpha_q,Q,\alpha_q,1,\alpha_q),\\ &(P,\alpha_q,Q,\alpha_q,1,\alpha_q), (P,\alpha_q,P,\alpha_q,1,\alpha_q),
(0,1,P + Q,P + \alpha_q,0,1), (P + Q, Q + \alpha_q,0,1,0,1),
\\&(1,0,1,0,1,0),  (\alpha_q P + PQ,\alpha_q P + \alpha_q Q,1,0,1,0), (1,0,\alpha_q Q + QP, \alpha_q P + \alpha_q Q,1,0)
\}
\end{align*}

Since all of them are rational double points, Lemma \ref{lemm-arith-gen} allows us to conclude that $p_a(S_q) = q + 1$ and $b_2(S_q) = 12q +10$.
\end{proof}

\begin{conj}
	Let $q$ be a prime power and let $S_q := \overline{\Sht}^{2, \tr}(\Gamma_0(N))$ over the field $F := \F_q(P,Q)$, where $N = (0) + (1) + (-1) + (\infty)$ if $q$ is odd and $N = (0) + (1) + (\alpha_q) + (\infty)$ if $q>2$ is even. Then
	\begin{enumerate}[\normalfont(1)]
		\item The trivial lattice of $S_q$ has rank $8q + 6$.
		\item If $q$ is odd, then $S_q$ has $4q + 8$ singular fibers over $\overbar{F}$: four fibers of type $\textnormal{I}_{2q + 2}$ at $0, 1, -1$ and $\infty$, and $4q + 4$ fibers of type $\textnormal{I}_1$.
		
	\end{enumerate}
\end{conj}

Finally, we conjecture that the generic fiber of $S_q$ has Mordell--Weil rank $4$, independently of $q$. This would be a consequence of the Tate conjecture, which predicts a relation between the second Betti number and the Picard number of a surface $S$ over $\F_q$. Fix a prime $\ell \neq p$, and let $H$ be the subgroup of NS$(S)$ generated by divisor classes defined over $\F_q$. Then the Tate conjecture predicts that the rank of $H$ equals the multiplicity of $q$ as an eigenvalue of the map induced by Frob$_q$ on $H^2(\overbar{S}, \Q_\ell)$. A consequence is that $\rho(\overbar{S})$, the geometric Picard number of $S$, equals the number of eigenvalues of Frob$_q^*$ on $H^2(\overbar{S}, \Q_\ell)$ of the form $q$ times a root of unity.
Recall that, by the Weil conjectures, the characteristic polynomial of Frob$_q^*$ acting on $H^2(\overbar{S}, \Q_\ell)$ has integer coefficients and factors over $\C$ as $\dsty \prod_{i = 1}^{b_2}(1 - \alpha_iT)$, with $|\alpha_i| = q$; hence $q^2/\alpha_i$ is an eigenvalue whenever $\alpha_i$ is, and these two values are different except when $\alpha_i = \pm q$. Combining this fact with the Tate conjecture, it would follow that the difference $b_2(S) - \rho(\overbar{S})$ is always even. If $S$ is an elliptic surface then $b_2$ is even, so this would imply that $\rho(\overbar{S})$ is even as well. 
In particular, note that the Tate conjecture is known to hold true for K3 elliptic surfaces with a section \cite{ArSD}, so  if $S$ is such a surface, then
$\rho(\overbar{S})$ is known to be even.

If $S = \overline{\Sht}^{2, 0}(\Gamma_0(N))$, then Drinfeld's work implies that $H^2(\overbar{S}, \Q_\ell)$ decomposes as a direct sum of a subspace coming from cusp forms, of dimension $4$ times the number of cusp forms of level $\Gamma_0(N)$, and an Eisenstein subspace on which Frob$_q^*$ acts with eigenvalues of the form $q$ times a root of unity.
Combining this fact with the Tate conjecture, we conclude that, conjecturally,
\[ b_2(S) = \rho(\overbar{S}) + 4\#\{ \text{cusp forms of level }\Gamma_0(N)\}. \] 
If moreover $S$ is an elliptic surface with generic fiber $E$ and trivial lattice $T$, this would mean
\[ b_2(S) = \textnormal{rank}(T) + \textnormal{rank}_{\text{MW}}(\overbar{E}) + 4\#\{ \text{cusp forms of level }\Gamma_0(N)\}, \] 
where $\overbar{E}$ denotes the base change of the generic fiber $E$ of $S_q$ to the algebraic closure.
On the other hand, it follows from \cite[Proposition 7.1]{DF} that there are $q$ cusp forms for GL$_2$ over $\F_q(T)$ with level $\Gamma_0(N)$. Hence, using our formula for $b_2(S_q)$ and our conjectured formula for rank$(T)$, we can solve for the Mordell--Weil rank of the generic fiber of $S_q$, obtaining
\[ \textnormal{rank}_{\text{MW}}(\overbar{E}) = (12q + 10) - (8q + 6) - 4q = 4. \]

\section{Future Work}\label{sec-fw}
When the curve $X$ is the projective line, experimental data strongly suggests that the spaces Sht$^{2, \tr}(\Gamma(N))$, Sht$^{2, \tr}(\Gamma_1(N))$, and Sht$^{2, \tr}(\Gamma_0(N))$ are surfaces of general type whenever the degree of $N$ is at least $2, 3$, and $5$, respectively. It should be possible to prove this by finding models for these surfaces that have at most rational double points as singularities.

We would like to find explicit equations for Sht$^{2, \tr}_{\PP^1}(\Gamma(N))$, Sht$^{2, \tr}_{\PP^1}(\Gamma_1(N))$, and Sht$^{2, \tr}_{\PP^1}(\Gamma_0(N))$ in the cases where the support of $N$ contains a point of degree higher than one, and to study how some invariants of these surfaces, like the arithmetic genus, behave under finite extensions of the constant field $\F_q$.

%EXAMPLES
Another goal is to provide proofs for the conjectures made in Section \ref{sec-deg-4}. Lemma \ref{lemm-arith-gen} would confirm the values of the arithmetic genus and the Betti number $b_2$, provided that we could find models for $\overline{\Sht}^{2, \tr}(\Gamma_0(N))$ having at most du Val singularities. Alternatively, we could study how other singularities affect the arithmetic genus.
The claims about the other invariants seem harder to prove, since computing them requires to know the generic fiber of the elliptic fibration of $\overline{\Sht}^{2, \tr}(\Gamma_0(N))$.

We would also like to study the Tate conjecture for the surfaces $\overline{\Sht}^{2, \tr}(\Gamma_0(N))$ in Section \ref{sec-deg-4}. For example, if $N = (0) + (1) + (-1) + (\infty)$ and $q$ is odd, we would prove the conjecture if we were able to find $4$ generators of the free part of the Mordell--Weil group of the generic fiber of $\overline{\Sht}^{2, \tr}(\Gamma_0(N))$; since the expected rank does not depend on $q$ , these generators might be independent of $q$ as well. To approach the conjecture in the remaining cases we would first need a formula for the number of cusp forms for GL$_2$ over $\F_q$ of level $\Gamma_0(N)$, and then we would proceed analogously.

Another future direction is to compute $\text{Sht}^2_X$, as well as the stacks with level structures, in cases where the base curve $X$ is an elliptic curve over $\F_q$. A difficulty is that the theory of rank $2$ vector bundles over an elliptic curve, studied in \cite{Ati}, is much more involved than in the $\PP^1$ case (in particular, rank $2$ vector bundles over an elliptic curve can be indecomposable, even when the curve is defined over an algebraically closed field).
There are two possible ways to approach this problem: one is to study the spaces of global sections of vector bundles over an elliptic curve; another is to exploit the fact that if $\mathcal{E}$ is a vector bundle over an elliptic curve, then its projectivization $\PP(\mathcal{E})$ is an elliptic ruled surface, so one can use elementary transformations of $\PP(\mathcal{E})$ to study the corresponding elementary modifications of $\mathcal{E}$ (see \cite[Section V.2]{Har} and \cite[Chapters 2 and 5]{Fri}).

%MODULARITY
One motivation to study moduli spaces of shtukas of rank $2$ is that they play a role in the following modularity conjecture for elliptic curves defined over function fields:
\begin{conj}[Modularity Conjecture]\label{conj-mod}
	Let $X$ be a smooth, projective, geometrically irreducible curve with function field $K$.
	Let $E$ be an elliptic curve over $K$ with conductor $N$ and whose associated Galois representation is irreducible. Then there exists a special kind of correspondence between $E \times E$ and the compactified moduli space of shtukas $\overline{\Sht}^2(\Gamma_0(N))$.
\end{conj}

This conjecture is derived from Drinfeld's study of the cohomology of moduli spaces of shtukas, the K\"unneth formula, and the Tate conjecture, which is open in the relevant setting.
See \cite[Conjecture 6.1.1.]{thesis} for the precise statement. 

When $q=2$, $X = \PP^1_{\F_2}$ and $E$ is an elliptic curve over $\F_2(X)$, the conjecture has been verified in two cases by explicitly constructing the predicted correspondence: the case $N = 2(0) + (1) + (\infty)$ is due to Elkies and Weinstein \cite{Elk}, and the case $N = 3(0) + (\infty)$ is joint work of the author with Elkies and Weinstein \cite[Section 6.2]{thesis}.
We would like to find more evidence to support this modularity conjecture, possibly by constructing explicit correspondences in other cases. 

We would also like to verify Conjecture \ref{conj-mod} for examples in which the base curve $X$ is an elliptic curve, once we have developed explicit equations for the corresponding moduli space of shtukas. Unlike in the $\PP^1$ case, there may exist nonconstant elliptic curves over $X$ having good reduction everywhere and irreducible associated Galois representation. While these curves cannot admit a surjection from the Drinfeld modular curve $X_0(N)$, Conjecture \ref{conj-mod} has no hypothesis about the reduction type of $E$ at any place of $K$, so we expect that these curves will be modular in this new sense.
We would be specially interested in verifying the conjecture for examples of this kind.

%BIBLIOGRAPHY
\bibliographystyle{siam}
 
\end{document}